\documentclass[12pt]{article}
\usepackage{amsfonts}
\usepackage{amssymb}
\usepackage{mathtools}
\usepackage{verbatim}
\newtheorem{theorem}{Theorem}
\newtheorem{theorem*}{Theorem}
\newtheorem{corollary}{Corollary}
\newtheorem{proposition}{Proposition}
\newtheorem{definition}{Definition}
\newtheorem{definition*}{Definition}
\newtheorem{remark}{Remark}

\begin{document}


\begin{center}
{{\Large \bf On  packing of  Minkowski balls. I
} \\
}
\end{center}

\begin{center}
{\bf Nikolaj M. Glazunov } \end{center}

\begin{center}
{\rm Glushkov Institute of Cybernetics NASU, Kiev, } \\
{\rm Institute of Mathematics and Informatics Bulgarian Academy of Sciences }\\
{\rm  Email:} {\it glanm@yahoo.com }
\end{center} 

\bigskip

{\bf Abstract.}

We investigate lattice packings of  Minkowski balls. By the results of the proof of Minkowski conjecture about the critical determinant 
we devide Minkowski balls on 3 classes: Minkowski balls, Davis balls and  Chebyshev-Cohn balls.
We investigate lattice packings of  these balls on  planes with varieng Minkowski metric and  search among these packings the optimal packings. In this paper we prove that the optimal lattice packing of the Minkowski, Davis, and Chebyshev-Cohn balls is realized with respect to the sublattices of index two of the critical lattices of
corresponding  balls.\\

{\bf Keywords:} lattice packing, Minkowski ball, Minkowski metric, critical lattice, optimal lattice packing.\\

{\bf 2020 Mathematics Subject Classification:}  11H06, 52C05 \\

{\bf Thanks.}{ The author is deeply grateful to the Bulgarian Academy of Sciences, the Institute of Mathematics and Informatics of the Bulgarian Academy of Sciences, Professor P. Boyvalenkov for their support.
The author was supported by Simons grant 992227.} \\

\section{Introduction}

A system of equal  balls in $n$-dimensional space is said to form a packing, if no two balls of the 
system have any inner points in common.

Lately the remarkable results in resolving the problem of optimal packing of balls in $8$ and 
$24$-dimensional real Euclidean spaces have been obtained \cite{via,ckmrv}.

In this research,  acknowledged by the Fields medal, optimal packings are constructed on lattices.

Our considerations connect with Minkowski conjecture~\cite{Mi:DA,M:LP,D:NC,Co:MC,W:MC,Ma:AC1,GM:MM,GM:P2,GGM:PM}  and use results of its proof.
Corresponding results and conjectures is most simply stated in terms of geometric lattices and critical
 lattices~\cite{Mi:GZ,Mi:DA,Cassels,lek}. The last are important partial case of geometric lattices.
We investigate lattice packings of  Minkowski balls, Davis balls  and Chebyshev-Cohn balls on  planes with varying Minkowski metric and  search among these packings the optimal packings.
The packing problem is studied on  classes of lattices 
related to the problem of the theory of Diophantine approximations 
considered by H. Minkowski \cite{Mi:DA} for the case of the plane.
For some other selected  problems and results of the theory Diophantine approximations see for instance \cite{AndersenDuke} and references therein.
The naming of balls is connecting with results of investigating of Minkowski conjecture on critical determinant and its justification is given below.
 In this paper we prove that the optimal lattice packing of the Minkowski, Davis, and Chebyshev-Cohn balls is realized with respect to the sublattices of index two of the critical lattices of
corresponding  balls.

\section{Minkowski conjecture, its proof and  Minkowski balls}

Let
 
$$ |\alpha x + \beta y|^p + |\gamma x + \delta y|^p \leq c \cdot 
 |\det(\alpha \delta - \beta \gamma)|^{p/2}, $$ 

be a diophantine inequality defined for a given real $ p >1 $;
hear $\alpha, \beta, \gamma, \delta$ are real numbers with 
$ \alpha \delta - \beta \gamma \neq 0 .$ 

H. Minkowski in his monograph~\cite{Mi:DA} raise the question
about minimum constant $c$ 
such that the inequality has integer 
solution other than origin. 
Minkowski with the help of his theorem on convex body has found a sufficient condition for the solvability of Diophantine inequalities in integers not both zero:
$$ c = {\kappa_p}^p,  \kappa_{p} = \frac{{\Gamma(1 + \frac{2}{p})}^{1/2}}{{\Gamma(1 + \frac{1}{p})}}.$$
 But this result is not optimal, and Minkowski also raised the issue of not improving constant $c$.
For this purpose Minkowski has proposed to use the critical determinant.

Recall the definitions~\cite{Cassels}.

Let $\mathcal R $ be a set and $\Lambda $ be a  lattice with base 
\[
\{a_1, \ldots ,a_n \}
\] 
in ${\mathbb R}^n.$ \\

A lattice $\Lambda $
is {\it admissible} for body $\mathcal R $ 
($ {\mathcal R}-${\it admissible})
if ${\mathcal R} \bigcap \Lambda = \emptyset $ or $0.$\\
 
 Let $\overline {\mathcal R}$ be the closure of $ {\mathcal R}$.
A lattice $\Lambda $
is {\it strictly admissible} for  $\overline {\mathcal R}$  
( $\overline {\mathcal R}-${\it strictly  admissible})
if $\overline {\mathcal R} \bigcap \Lambda = \emptyset $ or $0.$\\

Let
\[
  d(\Lambda) = |\det (a_1, \ldots ,a_n )|
 \] 
   be the determinant of 
$\Lambda.$ \\

The infimum
$\Delta(\mathcal R) $ of determinants of all lattices admissible for
$\mathcal R $ is called {\em the critical determinant} 
of $\mathcal R; $
if there is no $\mathcal R-$admissible lattices then puts
$\Delta(\mathcal R) = \infty. $ \\

 A lattice 
$\Lambda $ is {\em critical}
if $ d(\Lambda) = \Delta(\mathcal R).$ \\

 For the given 2-dimension
region $ D_p \subset {\mathbb R}^2 = (x,y), \ p \ge 1 $ :

$$ |x|^p + |y|^p < 1 , $$
let $\Delta(D_p) $ 
be the critical determinant of the region.

All determinants of admissible lattices of this domain that have three pairs of points on the boundary of this domain are parametrized
 by the Minkowski-Cohn moduli space of the form
 \begin{equation}
 \Delta(p,\sigma) = (\tau + \sigma)(1 + \tau^{p})^{-\frac{1}{p}}
  (1 + \sigma^p)^{-\frac{1}{p}}, 
 \end{equation}
in the domain
 $$ {\mathcal M}: \; \infty > p > 1, \; 1 \leq \sigma \leq \sigma_{p} =
 (2^p - 1)^{\frac{1}{p}}, $$
of the $ \{p,\sigma\} $-plane, where $\sigma$ is some real parameter \cite{Co:MC,GM:MM,GM:P2,GGM:PM,Gl4}. 

 In notations
\cite{GGM:PM} 
next result have proved:

\begin{theorem}
\cite{GGM:PM}
\label{ggm}
$$\Delta(D_p) = \left\{
                   \begin{array}{lc}
    \Delta(p,1), \; 1 \le p \le 2, \; p \ge p_{0},\\
    \Delta(p,\sigma_p), \;  2 \le p \le p_{0};\\
                     \end{array}
                       \right.
                           $$
here $p_{0}$ is a real number that is defined unique by conditions
$\Delta(p_{0},\sigma_p) = \Delta(p_{0},1),  \;
2,57 < p_{0}  < 2,58, \;  p_{0} \approx 2,5725$.
\end{theorem}
\begin{remark}
We will call $p_{0}$ the Davis constant.
\end{remark}
\begin{corollary}
$$ {\kappa_p} = {\Delta(D_p)}^{-\frac{p}{2}}.  $$
\end{corollary}

From Theorem (\ref{ggm}) in notations \cite{GGM:PM,Gl4}  we deduce the following corollary:  
\begin{corollary} (\cite{GGM:PM,Gl4})
\[
 {\Delta^{(0)}_p} = \Delta(p, {\sigma_p}) =  \frac{1}{2}{\sigma}_{p},   \; {\sigma}_{p} = (2^p - 1)^{1/p},
\]
\[
 {\Delta^{(1)}_p}  = \Delta(p,1) = 4^{-\frac{1}{p}}\frac{1 +\tau_p }{1 - \tau_p},   \;  2(1 - \tau_p)^p = 1 + \tau_p^p,  \;  0 \le \tau_p < 1. 
\]

  For their critical lattices respectively  $\Lambda_{p}^{(0)},\; \Lambda_{p}^{(1)}$ next conditions satisfy:   $\Lambda_{p}^{(0)}$ and 
 $\Lambda_{p}^{(1)}$  are  two $D_p$-admissible lattices each of which contains
three pairs of points on the boundary of $D_p$  with the
property that 
\[
(1,0) \in \Lambda_{p}^{(0)},\; (-2^{-1/p},2^{-1/p})
\in \Lambda_{p}^{(1)},
\]
 (under these conditions the lattices are
uniquely defined).
\end{corollary}

\section{Minkowski balls and the density  of the packing of $2$-dimensional Minkowski  balls}

Let $p_{0} \in {\mathbb R}$ be the Davis constant such that $2,57 < p_{0}  < 2,58$. \\

We consider balls of the form
\[
   D_p: \;  |x|^p + |y|^p \le 1, \; p \ge 1,
\]
and call balls with the conditions
\[
|x|^p + |y|^p \le 1,\;  2 > p > 1, 
\]
 the  Minkowski balls in two dimension and 
correspondingly the circles with  the condition
\[ 
|x|^p + |y|^p = 1, \; 2 > p > 1 
\]
the  Minkowski circles in two dimension.   
  
  Limiting Minkowski circle  in two dimension:    $|x| + |y| = 1$   \\

Davis balls  in two dimension: 
        $|x|^p + |y|^p \le 1$  for $p_{0} > p \ge 2$ \\

Davis circles in two dimension: 
    $|x|^p + |y|^p = 1$  for $p_{0} > p \ge 2$ \\

Chebyshev-Cohn balls in two dimension:   $|x|^p + |y|^p \le 1$  for $ p \ge p_{0}$ \\

Chebyshev-Cohn circles in two dimension: $|x|^p + |y|^p = 1$  for $ p \ge p_{0}$ \\

 Limiting Chebyshev-(Cohn) ball in two dimension: 
\[
||x, y||_{\infty} = \max(|x|,|y|).
 \]
 
Recall the definition of a packing lattice \cite{via,ckmrv,Cassels,lek}.
We will give it for $n$-dimensional Minkowski balls $D_p^n$ in ${\mathbb R}^n$.

\begin{definition}
 Let $\Lambda$ be a full lattice in ${\mathbb R}^n$ and $a \in {\mathbb R}^n$. 
 In the case if it is occurs that no two  of balls $\{D_p^n + b, b \in \Lambda + a \}$ have inner points in common, the collection of balls 
   $\{D_p^n + b, b \in \Lambda + a \}$ is called a $(D_p^n,  \Lambda)$-packing, and $\Lambda$ is called a packing lattice of $D_p^n$.
\end{definition}

Recall also that if $\alpha \in {\mathbb R}$ and $D_p^n$ is a ball than $\alpha D_p^n$ is the set of points $\alpha x, x \in D_p^n$.
 
In some cases we will consider interiors of balls $D_p =  D_p^2$ (open balls) which we will denoted as $ID_p$.

From the considerations of Minkowski and other authors \cite{Mi:DA,Mi:GZ,Cassels,lek,Gl4}, the following statements can be deduced 
(for the sake of completeness, we present the proof of Proposition \ref{p1}).\\

Denote by $V(D_p)$ the volume (area) of $D_p$.

\begin{proposition} \cite{Cassels,lek}.
\label{p1}
A lattice $\Lambda$ is a packing lattice of $D_p$ if and only if  it is admissible lattice for $2 D_p$.
\end{proposition}
{\bf Proof}  (contrary proof). First, note that one can take an open ball $ID_p$ and use the notion of strict admissibility.
 Suppose that $\Lambda$ is not strictly admissible for $2 ID_p$. Then $2 ID_p$ contains a point $a \ne 0$ of $\Lambda$.
Then the two balls $ID_p$ and $ID_p +a$ contain the point $\frac{1}{2} a$  in common. So $\Lambda$ is not a packing
 lattice of $D_p$.\\
 Suppose now that $\Lambda$ is not a packing  lattice of $D_p$. Then there exist two distinct points $b_1, b_2 \in \Lambda$ and a point  $c$ such that $c \in ID_p + b_1 $ and $c \in ID_p + b_2 $. 
 Hence there are points $a_1, a_2 \in ID_p$ such that $c = a_1 + b_1 = a_2 + b_2$. So $b_1 - b_2 = a_2 - a_1 \in ID_p$, whereas 
 $b_1 - b_2 \ne 0$ and $b_1 - b_2 \in \Lambda$. Therefore $\Lambda$ is not (strictly) admissible lattice.

\begin{proposition} \cite{Cassels,lek}.
\label{p3}
 The dencity of a $(D_p, \Lambda)$-packing is equal to $V(D_p)/d(\Lambda)$ and it is maximal if $\Lambda$ is critical for  $2 D_p$.
\end{proposition}

  

\section{On packing  Minkowski balls, Davis balls and  Chebyshev-Cohn balls on the plane}

Let us consider possible optimal lattice packings of these balls and their connection with critical lattices.

At first give lattices of trivial optimal lattice packings for the limiting  (asymptotic) cases at the points $p = 1$ and $\infty$ "infinity" (the latter corresponds to the classical Chebyshev balls) and as the introductory example the optimal packing of two-dimensional unit balls.

\begin{proposition}
\label{lp2}
The lattice 
\[
  \Lambda_{1}^{(1)} =  \{(\frac{1}{2}, \frac{1}{2}), (0, 1)\}
\]
is the critica lattice for $D_1$. 
Limiting case of Minkowski balls  for $ p = 1$ gives the optimal sublattice of index two of the lattice $\Lambda_{1}^{(1) }$-lattice  packing with the density $1$.
The centers of the Minkowski balls in this case are  at the vertices of the sublattice of index two of the lattice $\Lambda_{1}^{(1) }$.
\end{proposition}
{\bf Proof.} Recall that a critical lattice for $D_1$ is a lattice $\Lambda$ which is $D_1$-admissible and which has determinant 
$d(\Lambda) = \Delta(D_1)$. The lattice $\Lambda_{1}^{(1)}$ is $D_1$-admissible.   We have $\Delta(D_1) = \frac{1}{2}$ and $d(\Lambda_{1}^{(1)}) = \frac{1}{2}$. Minkowski balls  for $ p = 1$ are congruent squares. Hance we have the optimal 
the sublattice of index two of the lattice $\Lambda_{1}^{(1) }$  packing of the squares with the density $1$.

  \begin{proposition}
\label{lp4}
  The lattice 
\[
  \Lambda_{\infty}^{(1)} =  \{(1, 1), (0, 1)\}
 \]
is the critica lattice for $D_{\infty}$. 
Limiting case of Minkowski balls  for $ p = \infty$ gives  the optimal  of the  density $1$ packing 
 with respect to the sublattice of index two of the critical lattice $ \Lambda_{\infty}^{(1)}$.
 The centers of the Minkowski balls in this case are at the vertices of  the sublattice of index two of the lattice $\Lambda_{\infty}^{(1)}$.
  \end{proposition}
{\bf Proof.}    The lattice $\Lambda_{\infty}^{(1)}$ is $D_{\infty}$-admissible.   We have $\Delta(D_{\infty}) = 1$ and
 $d(\Lambda_{\infty}^{(1)} = 1$. Minkowski balls  for $ p = \infty$ are congruent squares. Hance we have the optimal 
the sublattice of index two of the lattice $\Lambda_{\infty}^{(1)}$  packing of the squares with the density $1$.\\

\begin{proposition}
\label{lp3}
The lattice 
\[
\Lambda_{2}^{(0)} =  \{(1, 0), (\frac{1}{2}, \frac{\sqrt{3}}{2})\}
\]
is the critica lattice for $D_2$. 
The lattice 
packing of Davis balls  for $p = 2$ gives the optimal  of the  density $\approx 0.91$ packing 
 with respect to the sublattice of index two of the critical lattice $\Lambda_{2}^{(0)}$.
 The centers of the Minkowski balls in this case are at the vertices of  the sublattice of index two of the lattice $\Lambda_{2}^{(0)}$.
\end{proposition}
{\bf Proof.} As in Proposition (\ref{lp2}) a critical lattice for $D_2$ is a lattice $\Lambda$ which is $D_2$-admissible and which has determinant 
$d(\Lambda) = \Delta(D_2)$. The lattice $\Lambda_{2}^{(0)}$ is $D_2$-admissible.  We have $\Delta(D_2) = \frac{\sqrt{3}}{2}$ and $d(\Lambda_{2}^{(0)}) = \frac{\sqrt{3}}{2}$. So sublattice of index two of the lattice $\Lambda_{2}^{(0)} $ is 
 the hexagonal lattice.  Next, we use the following classical results \cite{FegesToth,Thue}: the optimal  sphere packing of dimension 2  is the  hexagonal lattice (honeycomb) packing with the density $\approx 0.91$.
  
\begin{proposition}
\label{p2}
If $\Lambda$ is the critical lattice of the Minkowski ball $D_p$ than the sublattice $\Lambda_2$  of index two of the critical lattice is the critical lattice of $2 D_p$.
 (Examples for $ n = 1, 2, \infty $ above). Here we give the proof of the Proposition \ref{p2}.
\end{proposition}
{\bf Proof.} Since the Minkowski ball $D_p$ is symmetric about the origin and convex, then $2D_p$ is convex and symmetric
 about the origin \cite{Cassels,lek}.
 
When parametrizing admissible lattices $\Lambda$ having three pairs of points on the boundary of the ball $D_p$,
the following parametrization is used \cite{Co:MC,GM:P2,GGM:PM,Gl4}:
 \begin{equation}
\label{eq2}
\Lambda =  \{((1 + \tau^{p})^{-\frac{1}{p}}, \tau(1 + \tau^{p})^{-\frac{1}{p}}), (-(1 + \sigma^p)^{-\frac{1}{p}}, \sigma(1 + \sigma^p)^{-\frac{1}{p}})\}
\end{equation}
 where 
 \[
 0 \le  \tau < \sigma , \; 0 \le \tau \le \tau_p.
  \]
   $\tau_{p}$ is defined by the
equation $ 2(1 - \tau_{p})^{p} = 1 + \tau_{p}^{p}, \; 0 \leq
\tau_{p} < 1. $ 
\[
  1 \le\sigma \le \sigma_p, \; \sigma_p = (2^p - 1)^{\frac{1}{p}}.
\]
Admissible lattices of the form (\ref{eq2}) for doubled Minkowski balls $2 D_p$ have a representation of the form 
 \begin{equation}
\label{eq3}
\Lambda_{2 D_p} =  \{2((1 + \tau^{p})^{-\frac{1}{p}}, 2\tau(1 + \tau^{p})^{-\frac{1}{p}}), (-2(1 + \sigma^p)^{-\frac{1}{p}}, 2\sigma(1 + \sigma^p)^{-\frac{1}{p}})\}
\end{equation}
  Hence the Minkowski-Cohn moduli space for these admissible lattices has the form
 \begin{equation}
 \Delta(p,\sigma)_{2 D_p} = 4 (\tau + \sigma)(1 + \tau^{p})^{-\frac{1}{p}}
  (1 + \sigma^p)^{-\frac{1}{p}}, 
 \end{equation}
in the same domain
 $$ {\mathcal M}: \; \infty > p > 1, \; 1 \leq \sigma \leq \sigma_{p} =
 (2^p - 1)^{\frac{1}{p}}, $$
  
Consequently, the critical determinants of doubled Minkowski balls have a representation of the form
\[
 {\Delta^{(0)}_p}(2D_p) = \Delta(p, {\sigma_p})_{2D_p} =  2\cdot {\sigma}_{p},   \; {\sigma}_{p} = (2^p - 1)^{1/p},
\]
\[
 {\Delta^{(1)}_p}(2D_p)  = \Delta(p,1))_{2D_p} = 4^{1 - \frac{1}{p}}\frac{1 +\tau_p }{1 - \tau_p},   \;  2(1 - \tau_p)^p = 1 + \tau_p^p,  \;  0 \le \tau_p < 1. 
\]
  And these are the determinants of the sublattices of index 2 of the critical lattices of the corresponding Minkowski balls.
  
   \begin{theorem}
The optimal lattice packing of the Minkowski, Davis, and Chebyshev-Cohn balls is realized with respect to the sublattices of index two of the critical lattices
  \[
  (1,0)\in\Lambda_{p}^{(0)},\; (-2^{-1/p},2^{-1/p}) \in \Lambda_{p}^{(1)}.
  \]
    \end{theorem}
{\bf  Proof}.
 By Proposition \ref{p2} the critical lattice of $2 D_p$ is  the sublattice  of index two of the critical lattice 
of Minkowski ball $D_p$ .
So it is the admissible lattice for $2 D_p$ and by Proposition \ref{p1} 
is packing lattice of $D_p$.
By Proposition \ref{p3} the corresponding lattice packing has maximal density and so is optimal.
  
 \begin{remark}  
This result concerns the packing of unit balls and spheres in complete normed (Banach) spaces of dimension 2.
   \end{remark}

\end{document}